\documentclass[10pt]{article}

\usepackage{amscd,amstex,amssymb}

\long\def\comment#1\endcomment{}

\newcommand{\proof}[1][Proof]{\noindent {\it #1.\ }}
\def\endproof{\hfill \ensuremath{\square}\par}

\newcommand{\defn}{\noindent {\bf Definition.\ }}
\newcommand{\rem}{\noindent {\bf Remark.\ }}
\newcommand{\caution}{\noindent {\bf Caution.\ }}

\newenvironment{lemma}[1][]{\trivlist
   \item[\hskip \labelsep{\bfseries Lemma#1.}]\itshape}
   {\endtrivlist}

\newenvironment{corr}[1][]{\trivlist
   \item[\hskip \labelsep{\bfseries Corollary#1.}]\itshape}
   {\endtrivlist}

\newenvironment{prop}[1][]{\trivlist
   \item[\hskip \labelsep{\bfseries Proposition#1.}]\itshape}
   {\endtrivlist}

\newenvironment{theorem}[1][]{\trivlist
   \item[\hskip \labelsep{\bfseries Theorem#1.}]\itshape}
   {\endtrivlist}

\newenvironment{declare}[1][]{\trivlist
   \item[\hskip \labelsep{\bfseries #1.}]\itshape}
   {\endtrivlist}

\makeatletter

\def\@seccntformat#1{\csname the#1\endcsname.\hskip 5pt}

\def\punkt{\refstepcounter{subsection}
           \noindent{\bf \thesubsection.\ }}

\@addtoreset{equation}{section}

\newcommand{\ps@verbit}{%
  \renewcommand{\@oddhead}{%
          \scriptsize
          {Subvarieties of generalized Kummer varieties}
          \hfil\tiny {D. Kaledin, M. Verbitsky, \ \ \ January 7, 1998}}
  \renewcommand{\@evenhead}{\@oddhead}
  \renewcommand{\@oddfoot}{\hfil\thepage\hfil}
  \renewcommand{\@evenfoot}{\@oddfoot}}
 
\makeatother

\newcommand{\C}{{\Bbb C}}
\newcommand{\CC}{{\cal C}}
\newcommand{\Y}{{\frak Y}}
\newcommand{\X}{{\frak X}}
\newcommand{\Z}{{\cal Z}}
\newcommand{\N}{{\cal N}}
\newcommand{\T}{{\cal T}}

\newcommand{\wt}{\widetilde}
\newcommand{\wh}{\widehat}
\newcommand{\calo}{{\cal O}}
\newcommand{\sS}{{\frak S}}
\newcommand{\cp}{\C P^1}
\def\arrow{\to}
\newcommand{\birato}{\dashrightarrow}
\newcommand{\h}{{\Bbb H}}

\renewcommand{\dim}{\operatorname{\sf dim}} 
\newcommand{\Hilb}{\operatorname{\sf Hilb}}
\newcommand{\Diag}{\operatorname{\sf Diag}}
\newcommand{\graph}{\operatorname{\sf graph}}

\newcommand{\idot}{{\:\raisebox{3pt}{\text{\circle*{1.5}}}}}

\begin{document}

\begin{center}
{\Large\bf
Trianalytic subvarieties\\[2mm] of generalized Kummer varieties} \\[2mm]
D. Kaledin, M. Verbitsky,\\[4mm]
{\tt kaledin@@balthi.dnttm.rssi.ru} and
{\tt verbit@@thelema.dnttm.rssi.ru, verbit@@ihes.fr}
\end{center}

\hfill

\hspace{0.2\linewidth}
\begin{minipage}[t]{0.7\linewidth}
Let $X$ be a hyperk\"ahler manifold. Trianalytic subvarieties of $X$
are subvarieties which are complex analytic with respect to all
complex structures induced by the hyperk\"ahler structure.  Given a
2-dimensional complex torus $T$, the Hilbert scheme $T^{[n]}$ classifying
zero-dimensional subschemes of $T$ admits a hyperk\"ahler
structure. A finite cover of $T^{[n]}$ is a product of $T$ and a
simply connected hyperk\"ahler manifold $K^{[n-1]}$, called
generalized Kummer variety. We show that for $T$ generic, the
corresponding generalized Kummer variety has no trianalytic
subvarieties. This implies that a generic deformation of the
generalized Kummer variety has no proper complex subvarieties.
\end{minipage}

\tableofcontents

\section{Introduction}
In \cite{Beau} Beauville has constructed two series of examples of
compact hyperk\"ahler manifolds. The first one consists of the Hilbert
schemes $\Hilb$ of points on a K3 surface. The second is the series
of the so-called generalized Kummer varieties, one in each even
complex dimension. A generalized Kummer variety $K^{[n]}$ is
a subvariety of the Hilbert scheme of $n+1$ points on a
$2$-dimensional complex torus $K$, defined as a fiber of the
Albanese map.

Let $M$ be a hyperk\"ahler manifold. The hyperk\"ahler structure 
induces a 2-dimensional sphere of complex structures on $M$,
called {\bf induced complex structures.} A closed subset
of $M$ is called {\bf trianalytic} if it is complex analytic with respect
to all induced complex structures.

In \cite{Ver} it was proved that the Hilbert scheme of points on a
generic (in the sense of 
\ref{mumford.tate}) K3 surface has no proper trianalytic
subvarieties. In this paper we prove the same result for the
generalized Kummer varieties. Our methods are an adaptation of those
of \cite{Ver}, greatly simplified by the presence of a group
structure on a complex torus. We hope that this paper may serve as
an introduction to the more difficult situation of the K3 surface
treated in \cite{Ver}.

Let $M$ be a compact hyperk\"ahler manifold.
For a generic induced complex structure on $M$,
all complex analytic subvarieties of $M$ are trianalytic.
This implies that a generic deformation of a generalized
Kummer variety has no complex subvarieties. 

\subsection{Idea of the proof}

Here we give a rough sketch of the proof
of our main result. The rest of this paper
is independent from this Subsection.

Consider a generalized Kummer variety $K^{[n]}$.
To prove that $K^{[n]}$ has no trianalytic subvarieties,
we use the deformation theory of trianalytic subvarieties,
developed in \cite{VerSubvar}. 

Let $M$ be a compact hyperk\"ahler manifold and $I$ an induced complex
structure. We denote by $(M, I)$ the manifold $M$,
considered as a K\"ahler manifold. 
The cohomology of $M$ is equipped with 
a natural action of the group $SU(2)$ (see, e. g. \cite{VerTrian}).
Let $X\subset (M, I)$ be a closed complex subvariety. 
In \cite{VerTrian} (see also Theorem \ref{_triana_subva_SU(2)_Theorem_}),
it was proven that $X$ is trianalytic if and only if 
its fundamental class $[X]$ is $SU(2)$-invariant.
Since the fundamental class is deformationally invariant,
this implies that all complex analytic deformations of $X$
are trianalytic. A slightly more evolved version of this
argument (\cite{VerSubvar}) implies that the deformation
space of $X\subset (M, I)$ is hyperk\"ahler, and
the union of all deformations of $X\subset M$ is a trianalytic
subvariety of $M$. If the subvariety $X\subset (M, I)$
has no complex analytic deformations, it is called
{\bf rigid}. If $M$ contains
a proper trianalytic subvariety $X$, then either
$M$ contains a proper rigid trianalytic subvariety
$X'$ (obtained as a union of all deformations of $X$),
or $M$ has a finite cover which is a product of
two hyperk\"ahler manifolds (see Proposition
\ref{rigid_Proposition_}).
A generalized Kummer variety $K^{[n]}$ is simply connected,
and has $\dim H^{2,0}(K^{[n]})=1$. Therefore,
it cannot have such a finite cover. We obtain that
$K^{[n]}$ has a proper rigid trianalytic subvariety
if $K^{[n]}$ has a proper trianalytic subvariety.

The rest of the argument does not use the hyperk\"ahler
structure of $K^{[n]}$, but only uses the canonical holomorphic symplectic
structure. It is well known that a hyperk\"ahler manifold
$M$ is equipped with a canonical holomorphic symplectic structure
(\cite{Besse}). A trianalytic subvariety $X\subset (M, I)$
is {\bf non-degenerately
symplectic}, that is, the restriction of a holomorphic symplectic
form from $(M, I)$ to $X$ is non-degenerate outside of singularietes.
Given a non-degenerately symplectic subvariety $X\subset K^{[n]}$,
we show that $X$ is never rigid. Thus, $X$ cannot be trianalytic.

A generalized Kummer variety is canonically embedded into
a Hilbert scheme which is a desingularization of the symetric
power of a torus.
Consider the corresponding map $K^{[n]} \arrow T^{(n+1)}$
from the generalized Kummer variety to the $n+1$-th symmetric power of a
torus. We prove that, for $X\subset K^{[n]}$ 
a non-degenerately symplectic subvariety,
the map $X \stackrel \pi \arrow \pi(X)$ is finite
over the generic point of $\pi(X)$. This is done
using the basic properties of Hilbert schemes
(Lemma \ref{finite_Lemma_}).

Consider the subvariety $Y:= \pi^{-1}(\pi(X))\subset K^{[n]}$.
Over a generic point of $\pi(X)$, the map
$\pi:\; Y \arrow \pi(X)$ is a locally trivial
fibration. Using the group structure on $T^{n+1}$,
we show that this fibration admits a canonical
flat connection with finite monodromy, i. e.,
a trivialization over a finite covering (Lemma \ref{product_Lemma_}).  
To simplify notations, we assume for the duration of the Introduction 
that this trivialization is defined globally over the generic part of
$\pi(X)$. Denote by $F$ the generic fiber of $\pi:\; Y\arrow \pi(X)$. 
Let $\pi(X)_0$ denote the generic part of
$\pi(X)$, and $Y_0:= \pi^{-1}(\pi(X)_0)$.
Then $Y_0 = F \times \pi(X)_0$. Denote by 
$\Gamma\subset F \times \pi(X)$
the closure of $X\cap Y_0$ in $F \times \pi(X)$
Then $\Gamma$ is a correspondence between
$F$ and $\pi(X)$ which is finite over generic
points of $\pi(X)$. Given that $X$ is irreducible, we obtain
that $\Gamma$ is also irreducible. Assume now that
the complex structure on the torus $T$ is Mumford-Tate generic
(\ref{mumford.tate}). Then, by Theorem \ref{_triana_subva_SU(2)_Theorem_},
all complex subvarieties $Z\subset T^{n+1}$ are trianalytic.
Trianalytic subvarieties are hyperk\"ahler in the neighbourhood
of every smooth point. Therefore, the dimension $\dim_\C Z$
is even. We obtain that the variety $\pi(X)$ has no
complex subvarieties of codimension 1. Consider the
natural projection $p:\; \Gamma \arrow F$. 
Let $\cal C$ be a family of divisors passing through
every point of $F$. Unless $p(F)$ is a single point,
for each point $x\in \Gamma$ there exists a divisor
$C_x\in \cal C$ such that $p^{-1}(C_x)$ is a subvariety
of codimension 1 in $\Gamma$ which passes through $x$.
By construction of $\Gamma$, the natural projection
$p':\; \Gamma \arrow \pi(X)$ is finite over generic
point of $\pi(X)$.
Taking $x\in \Gamma$ generic, we obtain that the projection 
$p'(p^{-1}(C_x))$
of the divisor $p^{-1}(C_x)$ to $\pi(X)$ has codimension 
1 in $\pi(X)$. This gives a contradiction. Therefore,
$p(\Gamma)$ is a single point. We obtain that
$X\subset Y$ is the closure of a trivial section of a trivial
fibration $\pi:\; Y_0 \arrow \pi(X)_0$.
Such a section is determined by the choice of $f\in F$.
Varying the choice of $f$, we obtain a deformation of $X$.
Since $X$ is rigid, the map
$\pi:\; \pi^{-1}(\pi(X))\arrow \pi(X)$ is generically
finite. The fiber $\pi^{-1}(\pi(X))$ is a product of punctual
Hilbert schemes (Lemma \ref{product_Lemma_}). Therefore,
$\pi^{-1}(\pi(X))$ is connected, and 
consists of a single point. We obtain that
$X$ is an irreducible component of the
subvariety $\pi^{-1}(\pi(X))\subset K^{[n]}$.
Therefore, $X$ is rigid if and only if
$\pi(X)$ is a rigid subvariety of 
$\pi(K^{[n]})\subset T^{(n+1)}$.

Consider the natural action of $T$ on $T^{(n+1)}$.
Let $X'\subset T^{(n+1)}$ be a subvariety obtained
as a union of $t(\pi(X))$, for all $t\in T$.
Clearly, if $\pi(X)$ is rigid in $\pi(K^{[n]})$,
then $X'$ is rigid in $T^{(n+1)}$. To prove
that $K^{[n]}$ has no trianalytic subvarieties
it remains to study rigid subvarieties in
$T^{(n+1)}$.

Let $\sigma:\; T^{n+1} \arrow T^{(n+1)}$ be the natural
quotient map. Consider the variety $D:= \sigma^{-1}(X')$.
Since $T^{n+1}$ acts on itself by holomorphic automorphisms,
the variety $D$ is never rigid. Denote by $t(D)$ a deformation 
of $D$, associated with $t\in T^{n+1}$. 
By diagonals of $T^{n+1}$ we understand subvarieties given by
equations of type $x_i=x_j$ (see \eqref{_diagonals_Equation_}). 
Unless $D$ is a diagonal, 
$\sigma(t(D))\subset T^{(n+1)}$ is a deformation of $X'$, 
for appropriate $t\in T^{n+1}$. Therefore, all rigid subvarieties
of $T^{(n+1)}$ are diagonals. Since 
the map $\pi:\; \pi^{-1}(\pi(X))\arrow \pi(X)$ is generically
finite, $\pi(X)$ does not lie in the union of diagonals
of $T^{(n+1)}$. Therefore, $X'$ is not a diagonal. We obtain that
$\pi(X)$ is not rigid in $T^{(n+1)}$, and 
hense $X$ is not rigid in $K^{[n]}$. This gives a rough idea
of the proof of our result.


\section{Definitions and the statement of the theorem}


\punkt Recall that a hyperk\"ahler manifold $M$ is a Riemannian
manifold equipped an action of the quaternion algebra $\h$ in its
tangent bundle such that this action is smooth and parallel with
respect to the Levi-Civita connection. For excellent introductions
to hyperk\"ahler manifolds, we refer the reader to \cite{Besse} and
\cite{HKLR}.

Let $M$ be a hyperk\"ahler manifold. Every quaternion $h \in \h$
with $h^2=-1$ induces an almost complex structure on $M$. It is
well-known that all these almost complex structures are
integrable. We call them {\em the induced complex structures}. The
set of induced complex structures is naturally identified with the
complex projective line $\cp$. For every $I \in \cp$, we denote by
$M_I$ the manifold $M$ equipped with the corresponding induced
complex structure.

\punkt Every hyperk\"ahler manifold $M$ with any induced complex
structure $M_I$ is canonically holomorphically symplectic. Therefore
if $M$ is compact, then $\dim H^{2,0}(M_I) \geq 1$. A simply connected 
compact hyperk\"ahler manifold $M$ with $\dim H^{2,0}(M_I) = 1$ is called
{\em simple}. By a theorem of Bogomolov \cite{Bog}, every compact
hyperk\"ahler manifold has a finite covering which 
splits into a product of several
simple hyperk\"ahler manifolds and a complex torus.

A compact hyperk\"ahler manifold of complex dimension $2$ is either
a complex torus $T$, or a K3 surface $M$. Of these, only the K3
surfaces are simple. For every $n > 1$, examples of simple compact
hyperk\"ahler manifolds of complex dimension $2n$ were constructed
by Beauville in \cite{Beau}.

\punkt In this paper we study one of the two classes of
hyperk\"ahler manifolds introduced by Beauville, namely, the
so-called {\em generalized Kummer varieties}.  For the convenience
of the reader, we reproduce here their definition and main
properties.

Let $T$ be a complex torus of dimension $k$.  Consider the Hilbert
scheme $T^{[n+1]}$ of $n+1$ points on $T$. This is a complex variety
of dimension $k(n+1)$. The commutative group structure on the torus
$T$ defines a summation map $\Sigma:T^{n+1} \to T$, which induces a
summation map $\Sigma:T^{[n+1]} \to T$.

\defn The {\em generalized Kummer variety} $K^{[n]}$ associated to
the torus $T$ is the preimage $\Sigma^{-1}(0) \subset T^{[n+1]}$ of
the zero $0 \in T$ of the group structure on the torus $T$.

\punkt Assume that the complex torus $T$ is $2$-dimensional. In this
case the Hilbert scheme $T^{[n+1]}$ is smooth. The Kummer variety
$K^{[n]}$ associated to $T$ is also smooth. Moreover, it is
simply-connected, and $\dim H^{2,0}(K^{[n]}) = 1$.

Assume further that the torus $T$ is equipped with a hyperk\"ahler
structure. The holomorphic $2$-form associated to the hyperk\"ahler
structure on the torus $T$ defines a canonical non-degenerate
holomorphic $2$-form on the Hilbert scheme $T^{[n+1]}$. This form
gives by restriction a holomorphic symplectic form $\Omega$ on the
Kummer variety $K^{[n]}$. Therefore the canonical bundle of the
complex manifold $K^{[n]}$ is trivial. By the Calabi-Yau Theorem
\cite{Yau}, every K\"ahler class $\alpha \in H^{1,1}(K^{[n]})$
contains a unique Ricci-flat metric. By \cite{Beau}, every one of
these metrics together with the form $\Omega$ defines a
hyperk\"ahler structure on the Kummer variety $K^{[n]}$. The Kummer
variety equipped with any of these hyperk\"ahler structures is a
simple compact hyperk\"ahler manifold.

\caution There is no canonical choice for a K\"ahler structure on
the manifold $K^{[n]}$. Therefore, unlike the holomorphic symplectic
form, the hyperk\"ahler structure on the Kummer variety $K^{[n]}$ is
not defined by the hyperk\"ahler structure on the torus $T$.

\punkt We now recall some general facts on hyperk\"ahler manifolds
introduced in \cite{VerAct}, \cite{VerTrian} and
\cite{VerDesing}. Let $M$ be a hyperk\"ahler manifold, and let $X
\subset M$ be a closed subset.

\defn (\cite{VerTrian}) The subset $X \subset M$ is called {\em
trianalytic} if it is analytic for every induced complex structure
$I$ on $M$. 

\punkt Recall that every analytic subset $X \subset Y$ of dimension
$k$ in a compact complex manifold $Y$ of dimension $n$ defines a
canonical homology class $[X] \in H_{2k}(Y,\C)$ called {\em the
fundamental class} of the subset $X$. Using the Poincare duality
isomorphism $H_{2k}(Y,\C) \cong H^{2n-2k}(Y,\C)$, we can consider the
fundamental class $[X]$ as an element of the cohomology group
$H^{2n-2k}(Y,\C)$. 

Assume that the hyperk\"ahler manifold $M$ is compact. The
$\h$-action in the tangent bundle to $M$ induces a canonical
$SU(2)$-action in de Rham algebra of the manifold $M$. By
\cite{VerAct}, this action commutes with the Laplacian and induces
therefore an $SU(2)$-action in the cohomology $H^\idot(M,\C)$.  The
following criterion for trianalyticity is proved in \cite{VerTrian}.

\begin{theorem}[ (Trianalyticity criterion)]
\label{_triana_subva_SU(2)_Theorem_}
Let $I$ be an induced complex structure on $M$, and let $N \subset
M_I$ be a closed analytic subvariety of the complex manifold $M_I$.
Let $[N] \in H^\idot(M,\C)$ be the fundamental class of the
subvariety $N$. Then $N$ is trianalytic if and only if the
cohomology class $[N]$ is $SU(2)$-invariant.
\end{theorem}


\punkt Trianalytic subvarieties in hyperk\"ahler manifolds have many
special properties. Of these, the most important to us will be the
following theorem proved in \cite{VerDesing}.

\begin{declare}[Desingularization Theorem]
Let $X \subset M$ be a trianalytic variety in a compact
hyperk\"ahler manifold $M$, and let $I$ be an induced complex
structure on $M$. 

The normalization $\wh{X} \to M_I$ of the complex-analytic subvariety $X
\subset M_I$ is smooth, and the canonical projection $\wh{X} \to M$
induces a hyperk\"ahler structure on the smooth manifold $\wh{X}$.
\end{declare}

\punkt \label{mumford.tate} 
The goal of this paper is to study trianalytic subvarieties in the
Kummer variety associated to a generic hyperk\"ahler torus of
complex dimension $2$. The notion of genericity appropriate for our
purposes is the following one, introduced in \cite{Ver}.

\defn Let $X$ be a compact hyperk\"ahler manifold. An induced
complex structure $I$ on $X$ is called {\em Mumford-Tate generic}
with respect to the hyperk\"ahler structure if for all $n > 0$,
every cohomology class
$$
\alpha \in \bigoplus\limits_p H^{p,p}(X_I^n)\bigcap
H^{2p}(X^n,{\Bbb Z})\subset H^\idot(X^n,\C)
$$
is invariant under the canonical $SU(2)$-action.

As proved in \cite{VerGen}, for every compact hyperk\"ahler manifold
$X$ all the induced complex structures on $X$ except for a countable
number are Mumford-Tate generic. If an induced complex structure $I$
on a hyperk\"ahler manifold $X$ is Mumford-Tate generic, then it is
obviously also Mumford-Tate generic on any power $X^l$ of the
manifold $X$. Moreover, by the Trianalyticity Criterion of
\cite{VerTrian} every complex-analytic subvariety $Y \subset X_I$ is
trianalytic and, in particular, even-dimensional.

\punkt We can now formulate our main result. 

\begin{theorem}\label{main}
Let $T$ be a $2$-dimensional complex torus equipped with a
hyperk\"ahler structure. Assume that the complex structure on $T$ is
Mumford-Tate generic, and consider a generalized Kummer variety
$K^{[n]}$ associated to $T$. 

For any hyperk\"ahler structure on $K^{[n]}$ compatible with the
canonical holomorphic $2$-form, every irreducible trianalytic
subvariety $X \subset K^{[n]}$ is either the whole $K^{[n]}$ or a
single point.
\end{theorem}

\noindent
The proof of Theorem~\ref{main} takes up the rest of this paper. 

\punkt We finish this section with the following corollary of
Theorem~\ref{main}. 

\begin{corr}
Assume that the Kummer variety $K^{[n]}$ is equipped with a complex
structure $I$ which is Mumford-Tate generic with respect to some
hyperk\"ahler structure on $K^{[n]}$ compatible with the canonical
holomorphic $2$-form. 

Then every irreducible analytic subvariety $X \subset K^{[n]}_I$ is
either the whole $K^{[n]}$ or a single point. 
\end{corr}

\proof Indeed, by Theorem \ref{_triana_subva_SU(2)_Theorem_},
 every analytic subvariety in
$K^{[n]}_I$ is trianalytic. 
\endproof

\caution No matter which hyperk\"ahler structure on the Kummer
variety $K^{[n]}$ we take, the standard complex structure on
$K^{[n]}$ which comes from the embedding $K^{[n]} \subset T^{[n+1]}$
into the Hilbert scheme of the torus $T$ is not Mumford-Tate
generic.

\section{Reduction to the case of rigid subvarieties}

\punkt In this section we give our first reduction of
Theorem~\ref{main}. Namely, call a subvariety $X \subset Y$ in a
complex variety $Y$ {\em rigid\/} if it admits no local
deformations. In this section we prove the following.

\begin{prop}\label{rigid_Proposition_}
Let $Y$ be a complex manifold of dimension $n$ equipped with some
hyperk\"ahler structure. Assume in addition that the manifold $Y$ is
simply connected and that $\dim H^{2,0}(Y) = 1$.

If the manifold $Y$ admits a subvariety $X \subset Y$ of dimension
$k$, $0 < k < n$, which is trianalytic with respect to the
hyperk\"ahler structure on $Y$, then it also admits a rigid
subvariety of dimension $m$, $k \leq m < n$, trianalytic with
respect to the hyperk\"ahler structure on $Y$.
\end{prop}

All the generalized Kummer varieties $K^{[n]}$ are simply connected
and have $\dim H^{2,0}(K^{[n]}) = 1$. Thus to prove
Theorem~\ref{main} it suffices to prove that for a generic torus $T$
the associated Kummer variety $K^{[n]}$ has no proper rigid
subvarieties trianalytic with respect to some hyperk\"ahler
structure.

\punkt Before we prove Proposition~\ref{rigid_Proposition_}, we need
to recall several facts about moduli spaces of subvarieties in a
complex manifold. Let $Y$ be a compact complex manifold, and let $a
\in H^\idot(Y,\C)$ be a cohomology class of the manifold $Y$. Douady
\cite{Douady} had constructed a moduli space $S(Y,a)$ of complex
subvarieties in $Y$ with fundamental class $a$. The Douady moduli
space is equipped with the family $\X \to S(Y,a)$ and the universal
map $\X \to Y$ which coincide near every point $s \in S(Y,a)$ with
the universal family provided by the local deformation theory. If
the manifold $Y$ is K\"ahler, the Douady moduli space $S(Y,a)$ is
compact (\cite{_Fujiki_}, \cite{Lieber}).

\punkt The proof of Proposition~\ref{rigid_Proposition_} is based on the
following general fact. 

\begin{prop}\label{subvar}
Let $X$ and $Y$ be compact hyperk\"ahler manifolds.  Assume given an
immersion $f:X \to Y$ which is an embedding on a dense open subset
and compatible with the hyperk\"ahler structure. Fix a complex
structure $I$ on $Y$. Let $[f(X)]$ be the fundamental class of the
analytic subset $f(X_I) \subset Y_I$, and let $S = S(Y_I,[f(X)])$ be
the Douady moduli space of complex subvarieties in $Y_I$ with
fundamental class $[f(X)]$.

Then the complex-analytic space $S$ is a compact smooth
hyperk\"ahler manifold. Moreover, the total space $\wt{X}$ of the
universal family $\wt{X} \to S$ is also a smooth hyperk\"ahler
manifold, and the canonical map $\wt{X} \to Y$ is an immersion
compatible with the hyperk\"ahler structure. Finally, the projection
$\wt{S} \to S$ carries a canonical flat holomorphic connection.
\end{prop}

\punkt The proof of Proposition~\ref{subvar} begins in \ref{start}
and takes up the rest of this section. However, first we deduce
Proposition~\ref{rigid_Proposition_} from Proposition~\ref{subvar}. Indeed,
assume given a trianalytic subvariety $X_0 \subset Y$ in a complex
variety $Y$ equipped with some hyperk\"ahler structure.  Consider
the normalization $X \to X_0$ of the analytic subvariety $X_0
\subset Y$. By the Desingularization Theorem of \cite{VerDesing} the
normalization $X \to X_0$ is a smooth hyperk\"ahler manifold, and
the map $X \to Y$ is a trianalytic immersion. Moreover, it is an
embedding outside of the preimage in $X$ of the subset of singular
points of $X_0$. Therefore Proposition~\ref{subvar} applies to $X
\to Y$. Hence the universal family $\wt{X} = X \times S$ of
deformations of $X \to Y$ is a smooth hyperk\"ahler manifold, and
the canonical map $f:\wt{X} \to Y$ is an immersion.

Since $\wt{X}$ is compact, the image $f(\wt{X}) \subset Y$ is a
trianalytic subvariety, and $\wt{X} \to Y$ is the normalization of
the subvariety $f(\wt{X}) \subset Y$. We claim that the subvariety
$f(\wt{X}) \subset Y$ is rigid. Indeed, by \cite{VerSubvar} every
deformation $\wt{X}' \to Y$ of $\wt{X} \to Y$ is a hyperk\"ahler
manifold isometric (hence isomorphic) to $\wt{X}$. In particular, we
have a fibration $\wt{X}' \cong \wt{X} \to S$, and for every point
$s \in S$ the fiber $X'_s \to Y$ of this fibration over the point
$s$ is a deformation of the corresponding fiber $X_s \to Y$ of the
family $\wt{X} \to S$. Consequently the family $\wt{X}' \to Y$ is a
family of deformations of $X \to Y$. Since the family $\wt{X} \to Y$
is universal, we must have $\wt{X}' = \wt{X}$. Therefore the
subvariety $f(\wt{X}) \subset Y$ is rigid.

To prove Proposition~\ref{rigid_Proposition_}, it remains to show that $\dim
\wt{X}$ is strictly less than $n$. Indeed, assume that $\dim \wt{X}
= n$, so that the immersion $f:\wt{X} \to Y$ is \'etale. By
Proposition~\ref{subvar} the fibration $f:\wt{X} \to S$ carries a
flat holomorphic connection. Since its fibers are compact, we can
take a finite cover $S' \to S$ such that the pullback $\wt{X}
\times_S S'$ splits into a product
$$
\wt{X} \times_S S' \cong X \times S'.
$$
Since the manifold $Y$ is by assumption simply connected, the
\'etale map $X \times S' \to Y$ is an isomorphism. But by assumption
both $X$ and $S'$ are compact hyperk\"ahler manifolds. Therefore
$\dim H^{2,0}(X) \geq 1$ and $\dim H^{2,0}(S') \geq 1$. By the
K\"unneth formula $\dim H^{2,0}(X \times S') \geq 2$, which
contradicts $\dim H^{2,0}(Y) = 1$.  
\endproof

\punkt \label{start} Proposition~\ref{subvar} is a simple corollary
of the results of \cite{VerSubvar}, where an analogous statement was
proved for not necessarily smooth trianalytic subvarieties $X
\subset Y$. For the convenience of the reader, we sketch here an
alternative proof using the twistor spaces.

Let $f:X \to Y$ be as in the statement of Proposition~\ref{subvar},
and let $\pi:\Y \to \cp$ be the twistor space of the hyperk\"ahler
manifold $Y$. The complex manifold $Y_I$ is embedded into $\Y$ as
the fiber over $I \in \cp$. Moreover, as a smooth manifold, the
twistor space $\Y$ is canonically isomorphic to the product $\cp
\times Y$. Therefore by K\"unneth formula $H^\idot(\Y,\C) =
H^\idot(Y,\C) \otimes H^\idot(\cp,\C)$. 

Let $[f(X)]$ be the fundamental class of the complex subvariety
$f(X_I) \subset \Y$, and consider the Douady space $S(\Y,[f(X)])$ of
subvarieties in $\Y$ with fundamental class $[f(X)]$. Since the
projection $\pi:\Y \to \cp$ is proper, for every analytic subvariety
$X_s \subset \Y$ the image $\pi(X_s) \subset \cp$ is also an
analytic subvariety. Therefore it is either the whole $\cp$, or a
union of several points. 
 
Every point $s \in S(\Y,[f(X)])$ in the Douady space corresponds to
a subvariety $X_s \subset \Y$. Since the variety $\Y$ is proper over
$\cp$, the subset 
$$
S_{gen}(\Y,[f(X)]) \subset S(\Y,[f(X)])
$$ 
of points $s \in S(\Y,[f(X)])$ such that $\pi(X_s) = \cp$ is open in
the Douady space $S(\Y,[f(X)])$. Denote by $\sS \subset
S(\Y,[f(X)])$ the subset of points $s \in S(\Y,[f(X)])$ such that
the image $\pi(X_s) \subset \cp$ of the corresponding subvariety
$X_s \subset \Y$ consists of a single point. The subset $\sS \subset
S(\Y,[f(X)])$ is obviously a union of connected components of the
complement $S(\Y,[f(X)]) \setminus S_{gen}(\Y,[f(X)])$. Therefore it
is closed in the Douady space $S(\Y,[f(X)])$, hence also a complex
variety. In order to prove that the Douady space $S$ is
hyperk\"ahler, we will identify $\sS$ with its twistor space.

The correspondence $s \mapsto \pi(X_s) \in \cp$ defines a
holomorphic map $\pi:\sS \to \cp$.  Moreover, let $J = \pi(s) =
\pi(X_s) \in \cp$, so that we have $X_s \subset Y_J \subset \Y$.  By
the K\"unneth formula, the fundamental class $[X_s] \in H^\idot(Y_J)
= H^\idot(Y)$ in the cohomology of the fiber $Y_J$ coincides with
the class of $f(X)$ in the cohomology of $Y$. Therefore the fiber
$S_J = \pi^{-1}(\pi(s))$ of the space $\sS$ over the point $J \in
\cp$ coincides with the Doaudy space of subvarieties of $Y_J$ with
fundamental class $[f(X)]$. This applies, in particular, to the case
$I = J$, so that $S_I$ coincides with the space $S$.

\punkt We first prove that the variety $\sS$ is smooth and that the
projection $\pi:\sS \to \cp$ is submersive at every point $s \in
\sS$.

Let $X_s \subset \Y$ be the subvariety corresponding to the point
$s$. We have $[X_s] = [f(X)]$. By the Trianalyticity Criterion of
\cite{VerTrian} this implies that the submanifold $X_s \subset Y_J$
is trianalytic. Let $f_s:\wt{X}_s \to Y_J$ be its normalization. By
the Desingularization Theorem of \cite{VerDesing} the map $f_s$
induces on $\wt{X}_s$ the structure of a smooth hyperk\"ahler
manifold. Since the same applies to all deformations of $X_s \subset
\Y$ as well, the local universal moduli space for deformations of
the subvariety $X_s \subset \Y$ coincides with the local deformation
space for the pair $\langle \wt{X}_s, f_s\rangle$ of the smooth
manifold $\wt{X}_s$ and the map $f_s:\wt{X}_s \to \Y$.

Recall that we have a canonical short exact sequence 
$$
\begin{CD}
0 @>>> \T\left(\wt{X}_s\right) @>>> f_s^*(\T(\Y)) @>>> \N(f_s) @>>> 0
\end{CD}
$$
of holomorphic vector bundles on $\wt{X}_s$, where $\T(\wt{X}_s)$
and $\T(\Y)$ are tangent bundles, and $\N(f_s)$ is by definition the
normal bundle to the map $f_s$. By general deformation theory, the
formal completion of the universal local moduli space for
deformations of $f_s:\wt{X}_s \to \Y$ is isomorphic to the formal
neighborhood of $0$ in the certain cone $\CC_s$ in the space
$H^0(X_s,\N({f_s}))$ of global sections of the normal bundle
$\N(f_s)$. This cone is defined by the vanishing of the so-called
Massey products (\cite{_Illusie_}, \cite{_Retakh_}). 
In order to prove that the space $\sS$ is smooth at
$s \in \sS$, it suffices to prove that the Massey products vanish
identically, so that the cone $\CC_s$ coincides with the whole space
$H^0(X_s,\N({f_s}))$.

To prove this, we split the normal bundle $\N(f_s)$ in two pieces in
the following way.  Since $X_s \subset Y_J$, we can consider the map
$f_s:\wt{X}_s \to \Y$ as a map $f_s:\wt{X}_s \to Y_J$ into the fiber $Y_J
\subset \Y$ over the point $J \in \cp$. Let
$\N^\perp(f_s)$ be the normal bundle to $f_s:\wt{X}_s \to Y_J$, and
let $\CC^\perp_s \subset H^0(\N^\perp(f_s))$ be the Massey cone
corresponding to deformations of $f_s:X_s \to Y_J$ {\em inside} the
fiber $Y_J$.

The normal bundle $\N(f_s)$ splits canonically
$$
\N(f_s) = f_s^*\pi^*T_J\left(\cp\right) \oplus \N^\perp(f_s)
$$
into the sum of two bundles. The first is the pullback to $X_s$ of
the normal bundle to the fiber $Y_J \subset \Y$, which is isomorphic
to the constant rank-$1$ bundle $\pi^*T_J(\cp)$ whose fiber is the
tangent space $T_J(\cp)$ to $\cp$ at the point $J$. The second is
the normal bundle $\N^\perp(f_s)$ to the map $f_s:X \to Y_J$. This
splitting, in turn, induces a splitting
\begin{align*}
H^0(\wt{X}_s,\N(f_s)) &= H^0(\wt{X}_s,f_s^*\pi^*T_J(\cp)) \oplus
H^0(\wt{X}_s,\N^\perp(f_s))\\
&= T_J(\cp) \oplus H^0(\wt{X}_s,\N^\perp(f_s)), 
\end{align*}
and the cone $\CC_s$ lies in the product $T_J(\cp) \times \CC^\perp_s$. 

Now, by the Trianalyticity Criterion of \cite{VerTrian} 
(Theorem \ref{_triana_subva_SU(2)_Theorem_}) the
subvariety $X_s$ and all its deformations are trianalytic. Therefore
the cone $\CC_s$ contains the direct product $T_J(\cp) \times
\CC^\perp_s$, hence coincides with it. This proves that the
projection $\pi:\sS \to \cp$ is submersive at $s \in \sS$.

In order to prove that the variety $\sS$ is smooth at $s \in \sS$, it
remains to prove that the Massey products on the normal bundle
$\N^\perp_s(X_s)$ vanish, so that the cone $\CC^\perp_s$ is the whole space
$H^0(X_s,\N^\perp(f_s))$. This follows from the splitting of the
canonical exact sequence 
\begin{equation}\label{exact}
\begin{CD}
0 @>>> \T(\wt{X}_s) @>>> f_s^*(\T(Y_J)) @>>> \N^\perp(f_s) @>>> 0
\end{CD}
\end{equation}
of holomorphic bundles on $X_s$. This splitting was established in
\cite{VerSubvar}. As explained in \cite{VerSubvar}, it follows
from the fact that all these bundles are hyperholomorphic.

\punkt Recall that a pseudo-Riemannian manifold is a manifold
equipped with a nowhere degenerate symmetric 2-form in the tangent
bundle, not necessarily positively defined. A pseudo-K\"ahler 
manifold is a complex manifold equipped with a 
pseudo-Riemannian structure $(\cdot,\cdot)$, such that the complex
structure operator $I$ is orthogonal, and the corresponding
skew-symmetric 2-form $(\cdot,I \cdot)$ is closed.
A pseudo-hyperk\"ahler manifold is a pseudo-Riemannian
manifold equipped with three complex structure,
satisfying quaternion relations, which is pseudo-K\"ahler
with respect to these complex structures. For every
pseudo-hyperk\"ahler manifold, one can define its twistor
space, exactly as one does it for a hyperk\"ahler manifold.

\punkt By the general theory of twistor spaces developed in
\cite{HKLR}, to prove that $\sS$ is a twistor space for a
pseudo-hyperk\"ahler structure $\cal H$ on $S$, 
it remains to prove the following.
\begin{enumerate}
\item There exists an antiholomorphic involution $\sS \to
\overline{\sS}$ compatible with the antipodal involution on $\cp$. 
\item For every point $s \in \sS$ there exists a real section
$\wt{s}:\cp \to \sS$ of the projection $\pi:\sS \to \cp$ passing
through $s \in \sS$. 
\item For every such section $\wt{s}:\cp \to \sS$ the normal bundle
$\N(\wt{s})$ is a sum of several copies of the bundle $\calo(1)$ on
$\cp$. 
\item There exists a relative $\calo(2)$-valued non-degenerate
holomorphic $2$-form on $\sS$ over $\cp$.
\end{enumerate}
Defined this way pseudo-hyperk\"ahler structure $\cal H$
is unique.

\punkt The involution required in \thetag{a} is induced by the canonical
involution on the twistor space $\Y$. Indeed, this involution sends
subvarieties into subvarieties and acts identically on the subspace
$H^\idot(Y_I,\C) \subset H^\idot(\Y,\C)$. Therefore it preserves the
fundamental class $[f(X)] \in H^\idot(Y_I,\C) \subset
H^\idot(\Y,\C)$. 

The claim \thetag{b} follows from the fact that every deformation
$X_s$ of the manifold $X$ is trianalytic in its fiber $Y_J \subset
\Y$ of the twistor projection $\pi:\Y \to \cp$. Indeed, to obtain
such a real section, it suffices to take the twistor space $f_s:\X_s
\subset \Y$ of the desingularization $\wt{X}_s$ of the trianalytic
submanifold $X_s$ and let $\wt{s}:\cp \to \sS$ map a point $J' \in
\cp$ into the point in $\sS$ corresponding to the subvariety
$f_s(\pi^{-1}(J) \cap X_s)$.

The normal bundle $\N(\wt{s})$ to this real section coincides with
the direct image $\pi_*(\N(f_s))$ of the normal bundle $\N(f_s)$ to
the map $f_s:\X_s \to \Y$. Therefore \thetag{c} follows from
\cite{VerBundles}. Finally, under the identification 
$$
\N\left(\wt{s}\right) \cong \pi_*\left(\N(f_s)\right)
$$
the holomorphic $2$-form required in \thetag{d} is induced by the
canonical holomorphic $2$-form on the bundle $\N(f_s)$, and it is
non-degenerate by virtue of the splitting of the exact sequence
\eqref{exact}. 

To prove that $S$ is not only pseudo-hyperk\"ahler but
hyperk\"ahler, we need to check some positivity conditions, which is
easy.

\punkt To finish the proof of Proposition~\ref{subvar}, it remains
to prove that the universal family map $\wt{X} \to Y$ is an
immersion and to construct a flat holomorphic connection on the
fibration $\wt{X} \to S$. For this we refer the reader to
\cite{VerSubvar}, noting only that the splitting of the exact
sequence \eqref{exact} is crucial for both these facts.
\endproof

\section{Stratification by diagonals on a Hilbert scheme}

\punkt In order to proceed in our proof of Theorem~\ref{main}, we
need to recall several facts on the geometry of Hilbert schemes
$\Hilb$ of points on complex manifolds.

Let $M$ be a complex manifold of dimension $k$, and let $M^{(n)} =
M^{n} / \Sigma_{n}$ be the $n$-th symmetric power of the manifold
$M$. By definition the Hilbert scheme $M^{[n]}$ of $0$-dimensional
subschemes in $M$ of length $n$ maps into the space $M^{(n)}$. The
map $M^{[n]} \to M^{(n)}$ is proper. Its fibers are isomorphic to
products of punctual Hilbert schemes of dimension-$0$ subschemes of
$\C^k$ concentrated at $0$.

\punkt The variety $M^{(n)}$ is singular. However, it admits a
canonical stratification with non-singular strata. The strata of
this stratification are numbered by Young diagrams of length $n$,
that is, sequences of positive integers $k_1 \leq k_2 \leq ... \leq
k_l$ such that $\sum k_i = n$. The stratum $M^{(n)}_\Delta$
corresonding to a Young diagram $\Delta = \langle k_1 \leq\ldots\leq
k_l\rangle$ is by definition the subvariety in $M^{(n)}$ consisting
of orbits of points $\langle a_1,\ldots,a_n\rangle \in M^n$ such
that
\begin{align}\label{_diagonals_Equation_}
\begin{split}
a_1 &= a_2 = \ldots = a_{k_1}\\
a_{k_1+1} &= a_{k_1+2} = \ldots = a_{k_1+k_2}\\
&\ldots\\
a_{k_1+\ldots+k_{l-1}+1} &= a_{k_1+\ldots+k_{l-1}+2} = \ldots =
a_n,
\end{split}
\end{align}
and neither of $a_{k_1},a_{k_1+k_2},\ldots,a_n$ is equal to any
other. Every stratum $M^{(n)}_\Delta$ is smooth. It is isomorphic to
the quotient
$$
\left(M_1 \times M_2 \times \cdots \times M_l \setminus \Diag\right)/
\Sigma_\Delta, 
$$
where $M_1,\ldots,M_l$ are $l$ copies of the manifold $M$, $\Diag
\subset M_1 \times \cdots \times M_l$ is the subset of diagonals,
and $\Sigma_\Delta \subset \Sigma_l$ is the subgroup in the
symmetric group on $l$ letters consisting of transpositions which
fix the sequence $\Delta = \langle k_1,k_2,\ldots,k_l\rangle$.  We
will call this canonical stratification on the variety $M^{(n)}$
{\em the stratification by diagonals}.

\punkt The stratification by diagonals on the variety $M^{(n)}$
induces a stratification on the Hilbert scheme $M^{[n]}$. The strata
$M^{[n]}_\Delta$ are no longer necessarily smooth. The fiber of the
canonical proper map $M^{[n]}_\Delta \to M^{(n)}$ at a point
$\langle a_1,\ldots,a_n\rangle \in M^{(n)}$ is isomorphic to the
product of punctual Hilbert schemes of subschemes in $M$ of lengths
$k_1,\ldots,k_l$ concentrated at points $a_{k_1},\ldots,a_n \in M$.

Denote by $\eta:M^l \setminus \Diag \to M^{(n)}$ the canonical
Galois covering with the Galois group $\Sigma_\Delta$, and let
$$
\wt{M}^{[n]}_\Delta = M^{[n]}_\Delta \times_{M^{(n)}} \left(M^l
\setminus \Diag\right)
$$
be the pullback of the stratum $M^{[n]}_\Delta$ of the Hilbert
scheme $M^{[n]}$ with respect to this covering. 

\punkt The variety $\wt{M}^{[n]}$ admits a modular
interpretation. Indeed, it is the moduli space of pairs
$$
\begin{cases}
l \text{ different points }a_1,\ldots,a_l \in M \\
\text{For each }i, 1 \leq i \leq l, \text{ a subscheme
}\Z_i \subset M \text{ of length }k_i \text{ concentrated }\\
\text{at the point }a_i \in M.
\end{cases}
$$
This interpretation allows one to construct a canonical
compactification of the moduli space $\wt{M}^{[n]}_\Delta$.
Namely, it embeds as a dense open subset in the larger moduli space
of pairs
\begin{equation}\label{modular}
\begin{cases}
l \text{ points }a_1,\ldots,a_l \in M, \text{ not necessarily
different}\\ 
\text{For each }i, 1 \leq i \leq l, \text{ a subscheme
}\Z_i \subset M \text{ of length }k_i \text{ concentrated }\\
\text{at the point }a_i \in M.
\end{cases}
\end{equation}
This larger moduli space is compact. We denote it by
$\overline{M}^{[n]}_\Delta$.

\punkt Unfortunately, it seems that the canonical Galois covering
$\eta:\wt{M}^{[n]}_\Delta \to M^{[n]}_\Delta \subset M^{[n]}$ does
not extend to the compactification $\overline{M}^{[n]}_\Delta
\supset \wt{M}^{[n]}_\Delta$.

However, the map $\eta$ does extend to $\overline{M}^{[n]}_\Delta$
in a weaker sense. Recall that a {\em meromorphic map} $f:X
\birato Y$ from a complex variety $X$ to a complex variety $Y$ is by
definition an analytic subvariety $\Gamma \subset X \times Y$ such
that for an open dense subset $U \subset X$ the canonical projection
$$
\Gamma \bigcap \left(U \times Y\right) \to U
$$
is one-to-one. 

If the varieties $X$ and $Y$ are algebraic, then every algebraic map
$f:U \to Y$ from an open dense subset $U \subset X$ trivially
extends to a meromorphic map $X \birato Y$. However, for general
complex varieties this is no longer true. 

\begin{lemma}\label{proper}
The projection $\eta:\wt{M}^{[n]}_\Delta \to M^{[n]}_\Delta$ extends
to a meromorphic map $\overline{\eta}:\overline{M}^{[n]}_\Delta
\birato M^{[n]}$ from the compactification
$\overline{M}^{[n]}_\Delta$ into the Hilbert scheme $M^{[n]}$.
\end{lemma}

\proof Note that the lower arrow $\eta:\left(M^l \setminus \Diag\right)
\to M^{(n)}_\Delta$ in the Cartesian square
$$
\begin{CD}
\wt{M}^{[n]}_\Delta @>{\eta_0}>> M^{[n]}\\
@V{\pi}VV                          @VV{\pi}V\\
M^l \setminus \Diag   @>{\eta_0}>> M^{(n)},
\end{CD}
$$
does extend to a holomorphic map
$$
\eta:M^l \to M^{(n)}.
$$
Let $\overline{M}'_\Delta$ be the fibered product given by the
Cartesian square 
$$
\begin{CD}
\overline{M}'_\Delta @>{\eta}>> M^{[n]}\\
@VVV                            @VV{\pi}V\\
M^l                  @>{\eta}>> M^{[n]}.
\end{CD}
$$
Both varieties $\overline{M}'_\Delta$ and
$\overline{M}^{[n]}_\Delta$ are proper over the complex manifold
$M^l$. Moreover, the preimage of the open subset
$\left(M^l\setminus\Diag\right) \subset M^l$ in each of these
varieties is canonically isomorphic to the variety
$\wt{M}^{[n]}_\Delta$.

Since by construction we have a map $\eta:\overline{M}'_\Delta \to
M^{[n]}$, it suffices to prove that the embedding
$\iota:\wt{M}^{[n]}_\Delta \hookrightarrow \overline{M}'_\Delta$
extends to a meromorphic map $\iota:\overline{M}^{[n]}_\Delta \birato
\overline{M}'_\Delta$. In other words, we have to prove that
\begin{itemize}
\item the closure $\overline{\graph\iota}
\subset \overline{M}^{[n]}_\Delta \times \overline{M}'_\Delta$ of
the graph
$$
\graph\iota \subset \wt{M}^{[n]}_\Delta \times \overline{M}'_\Delta
\subset \overline{M}^{[n]}_\Delta \times \overline{M}'_\Delta
$$
of the embedding $\iota:\wt{M}^{[n]}_\Delta \hookrightarrow
\overline{M}'_\Delta$ is an analytic subvariety. 
\end{itemize}
The statement \thetag{$\bullet$} is obviously local on $M^l$,
therefore it suffices to prove it in a neighborhood of every point
$m \in M^l$. Moreover, it is non-trivial only over a neighborhood of
the subset of diagonals $\Diag \subset M^l$. By induction, it
suffices to prove it over the neighborhood of the smallest diagonal
$M \subset M^l$, or, in other words, over subsets of the form $U^l
\subset M^l$ for a neighborhood $U \subset M$ of every point $m \in
M$.

But all the varieties and maps in \thetag{$\bullet$} admit modular
interpretations as Hilbert schemes of points on $M$ with some
additional conditions. Therefore everything in \thetag{$\bullet$} is
functorial with respect to open embeddings. Consequently, for every
open subset $U \subset M$ the statement \thetag{$\bullet$} holds
over $U^l \subset M^l$ if and only if it hold for $U$ in place of
$M$.

It remains to notice that, since all complex manifolds of the same
dimension $k$ are locally isomorphic, the statement
\thetag{$\bullet$} holds for every one of them if and only if it
holds for an arbitrary single manifold $M$. However, for the
algebraic manifold $M = \C^k$ the statement \thetag{$\bullet$} holds
trivially. Therefore it must hold for every complex manifold $M$ of
dimension $k$.  
\endproof

\section{Stratification of the Kummer variety}

\punkt Everything said in the last section was valid for an
arbitrary complex manifold $M$. The main result of this section,
Lemma~\ref{product_Lemma_}, is specific for the case of a complex torus. It
is this result that makes the proof of Theorem~\ref{main} for the
Kummer varieties much simpler than the proof of the analogous
statement for K3 surfaces given in \cite{Ver}.

\punkt Let $T$ be a complex torus. Let $T^{[n+1]}$ be the Hilbert
scheme of $0$-di\-men\-si\-o\-nal subschemes of length $n+1$ in the
torus $T$, and let
$$
K^{[n]} \subset T^{[n+1]}
$$
be the associated generalized Kummer variety. 

The summation map $\Sigma:T^{n+1} \to T$ factors through the
symmetric power $T^{(n+1)} \to T$.  Denote by $K^{(n)} \subset
T^{(n+1)}$ the preimage $\Sigma^{-1}(0) \subset T^{(n+1)}$ of the
zero $0 \in T$. The canonical projection $\pi:T^{[n+1]} \to
T^{(n+1)}$ commutes with the summation map and defines therefore a
proper map $\pi:K^{[n]} \to K^{(n)}$.

\punkt The stratification by diagonals on the variety $T^{(n+1)}$
induces a stratification on the varieties $K^{(n)}$ and
$K^{[n]}$. For every Young diagram $\Delta$ of length $n+1$, denote
by
\begin{align*}
K^{(n)}_\Delta &= T^{(n+1)}_\Delta \bigcap K^{[n]} \subset T^{(n+1)}\\
K^{[n]}_\Delta &= T^{[n+1]}_\Delta \bigcap K^{[n]} \: \subset T^{[n+1]}
\end{align*}
the corresponding strata of this stratification. 

The canonical Galois covering $\eta:\left(T^l\setminus\Diag\right)
\to T^{(n+1)}_\Delta$ also commutes with the summation
map. Therefore it defines a Galois covering
$$
\eta:T^l_0\setminus\Diag \to K^{(n)}_\Delta, 
$$
where $T^l_0 \subset T^l$ is the kernel of the summation map
$\Sigma:T^l \to T$. Let 
$$
\wt{K}^{[n]} = \left(T^l_0\setminus\Diag\right)
\times_{K^{(n)}} K^{[n]} \subset \wt{T}^{[n+1]}_\Delta
$$
be the pullback variety. The summation map
$\Sigma:\wt{T}^{[n+1]}_\Delta \to T$ extends to the compactification
$\overline{T}^{[n+1]}_\Delta \supset \wt{T}^{[n+1]}_\Delta$, which
gives a compactification
$$
\wt{K}^{[n]}_\Delta \subset \overline{K}^{[n]}_\Delta =
\Sigma^{-1}(0) \subset \overline{T}^{[n+1]}_\Delta.
$$ 

\punkt Lemma~\ref{proper} immediately implies that the map
$\eta:\wt{K}^{[n]}_\Delta \to K^{[n]}_\Delta$ extends to a
meromorphic map $\eta:\overline{K}^{[n]}_\Delta \birato K^{[n]}$. We
will need the following corollary of this lemma.

\begin{corr}\label{image}
Let $X \subset \wt{K}^{[n]}_\Delta$ be an analytic subset. If the
closure $\overline{X} \subset \overline{K}^{[n]}_\Delta$ is an
analytic subset in the variety $\overline{K}^{[n]}_\Delta$, then the
closure
$$
\overline{\eta(X)} \subset K^{[n]}
$$
of the image $\eta(X) \subset K^{[n]}$ is also an analytic
subvariety. 
\end{corr}

\proof Indeed, since the meromorphic map
$\overline{\eta}:\overline{K}^{[n]}_\Delta \to K^{[n]}$ is given by
a proper correspondence, the direct image
$\overline{\eta}(\overline{X}) \subset K^{[n]}$ is an analytic
subvariety. The closure $\overline{\eta(X)} \subset K^{[n]}$ is a
union of irreducible components of the subvariety
$\overline{\eta}(\overline{X})$.  
\endproof

\punkt For any integer $l \geq 2$, denote by $F_l$ the punctual
Hilbert scheme of length $l$, that is, the moduli space of
$0$-dimensional subschemes in $\C^k$ of length $l$ concentrated at
$0$. For any Young diagram $\Delta = \langle k_1 \leq \cdots \leq
k_l\rangle$ let $F_\Delta = F_{k_1} \times \cdots \times F_{k_l}$ be
the product of $l$ such Hilbert schemes of lengths $k_1,\ldots,k_l$.
The main result of this section is the following.

\begin{lemma}\label{product_Lemma_} 
There exists a direct product decomposition
\begin{equation}\label{looking.for}
\overline{K}^{[n]}_\Delta \cong F_\Delta \times T^l_0. 
\end{equation}
\end{lemma}

\proof Recall that the variety $\overline{T}^{[n+1]}_\Delta$ is the
moduli space of data \eqref{modular}. The group $T^l$ acts on these
data by left translations, which induces a $T^l$-action on the
moduli space $\overline{T}^{[n+1]}_\Delta$. This action is free and
gives a decomposition
$$
\overline{T}^{[n+1]}_\Delta \cong F_\Delta \times T^l. 
$$
This decomposition is compatible with the summation map and induces
the desired decomposition \eqref{looking.for}.  
\endproof

\section{Subvarieties of a stratum $K^{[n]}_\Delta$}

\punkt Up to this point we did not need any facts on the torus $T$
except for its group structure. From now on, assume that the complex
torus $T$ is $2$-dimensional and equipped with a hyperk\"ahler
structure. Moreover, assume that the complex structure on $T$ is
Mumford-Tate generic with respect to this hyperk\"ahler structure in
the sense of \ref{mumford.tate}. In particular, there are no
analytic subvarieties in $T$ except for $T$ itself and unions of its
points.

Since $\dim T = 2$, the Hilbert scheme $T^{[n+1]}$ and the Kummer
variety $K^{[n]}$ are smooth. The hyperk\"ahler structure on the
torus $T$ induces a natural holomorphic symplectic structure on the
associated Kummer variety $K^{[n]}$. Fix once and for all a
hyperk\"ahler structure on $K^{[n]}$ compatible with the canonical
holomorphic symplectic form.

\punkt The hyperk\"ahler structure on the torus $T$ induces a
canonical hyperk\"ahler structure on the powers $T^l$ of $T$. The
first consequence of the genericity of the hyperk\"ahler structure
on the torus $T$ is the following.

\begin{lemma}\label{all.trian}
Every analytic subvariety $X \subset T^l$ is trianalytic.
\end{lemma}

\proof By the Trianalyticity Criterion of \cite{VerTrian} it
suffices to prove that the fundamental class $[X] \in
H^\idot(T^l,\C)$ is invariant under the canonical
$SU(2)$-action. But since $T$ is Mumford-Tate generic, every Hodge
cohomology class in $H^\idot(T^l,\C)$ is $SU(2)$-invariant.
\endproof

\punkt Applying the theory of trianalytic subvarieties developed in
\cite{VerSubvar}, we get the following stronger statement. 

\begin{lemma}\label{generic}
Let $X \subset T^l$ be an analytic subvariety, and denote by
$$
\pi:\C^{2l} \to T^l
$$
the universal covering map. 

Then every irreducible component of the variety $X$ is a complex
torus isogenous to a power of $T$. Moreover, the subvariety $X
\subset T^l$ is of the form $\pi(V)$, where $V \subset \C^{2l}$ is a
union of hyperplanes.
\end{lemma}

\proof Indeed, since the subvariety $X \subset T^l$ is trianalytic
by Lemma~\ref{all.trian}, then so is the subvariety $\pi^{-1}(X)
\subset \C^{2l}$. By \cite{VerSubvar}, Corollary 5.4, every trianalytic
subvariety in a hyperk\"ahler manifold is
totally geodesic.  A totally geodesic subvariety in
a flat manifold is a union of hyperplanes.  
\endproof

\punkt Let now $X \subset K^{[n]}$ be an analytic subvariety which
is trianalytic with respect to the chosen hyperk\"ahler structure on
the Kummer variety $K^{[n]}$. Say that the subvariety $X$ {\em
generically} lies in the stratum $K^{[n]}_\Delta \subset K^{[n]}$ if
\begin{enumerate}
\item $X$ lies in the closure of the stratum $K^{[n]}_\Delta \subset
K^{[n]}$, and 
\item the intersection $X \cap K^{[n]}_\Delta \subset X$ is a dense
open subset. 
\end{enumerate}
Every irreducible subvariety $X \subset K^{[n]}$ obviously lies
generically in one and only one stratum of the stratification by
diagonals.

\punkt In order to apply the genericity of the hyperk\"ahler
structure on $T$ to the study of the subvariety $X \subset K^{[n]}$,
we first prove the following.

\begin{lemma}\label{finite_Lemma_}
Assume that the trianalytic subvariety $X \subset K^{[n]}$ lies
generically in the stratum $K^{[n]}_\Delta$, and let $\pi:K^{[n]}
\to K^{(n)}$ be the canonical proper map. Then the induced
projection $\pi:X \to \pi(X)$ is finite and \'etale over an open
dense subset $U \subset \pi(X)$.
\end{lemma}

\proof It suffices to take $U \subset \pi(X) \cap K^{[n]}_\Delta$
and such that $\pi^{-1}(U) \subset X$ is smooth. Then the
hyperk\"ahler structure on $K^{[n]}$ induces a hyperk\"ahler
structure on $\pi^{-1}(U)$. Therefore the restriction
$\Omega\mid_{\pi^{-1}(U)}$ to $\pi^{-1}(U)$ of the canonical
holomorphic $2$-form $\Omega$ on $K^{[n]}$ associated to the
hyperk\"ahler structure must be non-degenerate. 

On the other hand, the canonical holomorphic $2$-form associated to
the hyperk\"ahler structure on the torus $T$ gives a holomorphic
$2$-form on the variety $T^l_0$. This form is invariant with respect
to the action of the group $\Sigma_\Delta$ on $T^l_0$. Therefore it
induces a holomorphic $2$-form $\wt{\Omega}$ on the smooth complex
manifold $K^{(n)}_\Delta = \left(T^l_0 \setminus
\Diag\right)\times\Sigma_\Delta$. By assumption the subvariety
$\pi^{-1}(U) \subset K^{[n]}_\Delta$ is smooth, therefore by
\cite{Ver}, Proposition 4.5, we have
$$
\Omega\mid_{\pi^{-1}(U)} = \pi^*\wt{\Omega}\mid_U. 
$$
Since this form is non-degenerate, the projection $\pi:U \to \pi(U)
\subset K^{(n)}_\Delta$ must be unramified. Since it is proper,
it must therefore, indeed, be finite and \'etale. 
\endproof

\punkt \label{wtU} 
Let now $X \subset K^{[n]}$ be an irreducible trianalytic subvariety
which lies generically in the stratum $K^{[n]}_\Delta \subset
K^{[n]}$. By the Desingularization Theorem of \cite{VerDesing} the
normalization $\wh{X}$ of the variety $X$ is smooth, and the
normalization map $\nu:\wh{X} \to X \hookrightarrow K^{[n]}$ 
is an immersion which induces
on $\wh{X}$ a hyperk\"ahler structure.

Choose a dense open subset $U \subset \pi(X)$ as in
Lemma~\ref{finite_Lemma_}, and let 
$$
\wh{U} = \nu^{-1}\left(X \bigcap \pi^{-1}(U)\right) \subset \wh{X}  
$$
be its preimage in $\wh{X}$. Moreover, consider the pullback $\wh{U}
\times_{K^{[n]}_\Delta} \wt{K}^{[n]}_\Delta$ of the variety $\wh{U}
\to K^{[n]}_\Delta$ with respect to the canonical Galois covering
$\eta:\wt{K}^{[n]}_\Delta \to K^{[n]}_\Delta$, and let $\wt{U}$ be
any one of its irreducible components. 

Recall that by assumption the subvariety $X \subset K^{[n]}$ is
irreducible. Therefore both its normalization $\wh{X}$ and the open
dense subset $\wh{U} \subset \wh{X}$ are irreducible, which implies
that $\eta(\wt{U}) \subset K^{[n]}_\Delta \subset K^{[n]}$ is dense
in the analytic subvariety $X \subset K^{[n]}$.

\punkt By definition the variety $\wt{U}$ is equipped with a
canonical map $\rho:\wt{U} \to \wt{K}^{[n]}_\Delta$. Under the
product decomposition \eqref{looking.for} the map $\rho$ decomposes
$\rho = \rho_T \times \rho_F$ into a product of a holomorphic map
$$
\rho_T:\wt{U} \to T^l_0 \setminus \Diag
$$
and a holomorphic map 
$$
\rho_F:\wt{U} \to F_\Delta.
$$
Moreover, by Lemma~\ref{finite_Lemma_} the first factor
$\rho_T:\wt{U} \to T^l_0 \setminus \Diag$ in this decomposition is
\'etale onto its image. The main result of this section is the
following.

\begin{prop}\label{point}
The canonical map $\rho_F:\wt{U} \to F_\Delta$ is a projection onto
a single point $o \in F_\Delta$.
\end{prop}

\proof Let $\overline{U} \subset \overline{K}^{[n]}_\Delta$ be the
closure of the subvariety $\wt{U} \subset K^{(n)}_\Delta$ in the
compactification $\overline{K}^{[n]}_\Delta$. This closure is an
analytic subvariety. Indeed, since the manifold $K^{[n]}$ is proper,
the preimage $\eta^{-1}(X) \subset \overline{K}^{[n]}_\Delta$ of the
subvariety $X \subset K^{[n]}$ under the meromorphic map
$\eta:\overline{K}^{[n]}_\Delta \birato K^{[n]}$ is an analytic
subvariety. Since the variety $\wt{U}$ is irreducible, its closure
$\overline{U} \subset \overline{K}^{[n]}_\Delta$ coincides with an
irreducible component of the subvariety $\eta^{-1}(X) \subset
\overline{K}^{[n]}_\Delta$.
 
Assume that the projection $\rho_F:\overline{U} \to F_\Delta$ does
not map $\wt{U} \subset \overline{U}$ to a single point. Since
$\wt{U}$ is irreducible, this implies that $\rho_F(\wt{U}) \subset
\rho_T(\overline{U}) \subset F_l$ is a variety of positive
dimension. But the punctual Hilbert scheme $F_\Delta$ is
projective. Thus we can take the preimage of the appropriate
hyperplane section in $F_l$ and obtain an analytic subvariety $D
\subset \overline{U}$ such that the intersection $D \cap \wt{U}
\subset \wt{U}$ is a non-trivial subvariety of codimension $1$.

Since the projection $\rho_T:\overline{U} \to T^l_0$ is proper, the
image $\rho_T(\overline{U}) \subset T^l_0$ is an analytic subvariety
in the torus $T^l_0$. The image $\rho_T(D) \subset
\rho_T(\overline{U})$ is also analytic. Moreover, the map
$\rho_T:\wt{U} \to \rho_T(\wt{U})$ is finite. Therefore the
intersection $\rho_T(D) \cap \rho_T(\wt{U}) \subset \rho_T(\wt{U})$
is a subvariety of codimension $1$. Since its closure
$$
\overline{D} = \overline{\rho_T(D) \cap \rho_T(\wt{U})} \subset
\rho_T(\overline{U})
$$ 
is a union of irreducible components of the analytic subvariety
$\rho_T(D) \subset \rho_T(\overline{U})$, it is also an analytic
subvariety, and its codimension is exactly $1$.

Now, by Lemma~\ref{all.trian} the subvariety $\rho(\overline{U})
\subset T^l_0$ is trianalytic. Therefore by Lemma~\ref{generic}
every one of its irreducible components is a complex torus isogenous
to a power of the torus $T$. Since the torus $T$ is Mumford-Tate
generic in the sense of \ref{mumford.tate}, all its analytic
subvarieties are even-dimensional, and the same is true for its
powers $T^l$. This contradicts the existense of $\overline{D}
\subset \rho(\overline{U})$ and proves the proposition.  \endproof

\section{Proof of the main theorem} 

\punkt We can now prove Theorem~\ref{main}. Let $X \subset K^{[n]}$
be an analytic subvariety in the Kummer variety $K^{[n]}$ which is
trianalytic with respect to the chosen hyperk\"ahler structure on
$K^{[n]}$. Assume that $\dim X > 0$. We have to show that $X =
K^{[n]}$.

By Proposition~\ref{rigid_Proposition_} we can assume that the subvariety $X
\subset K^{[n]}$ is rigid. Assume also, by induction, that every
trianalytic subvariety $Y \subset K^{[n]}$ with $\dim Y > \dim X$
coincides with the whole Kummer variety $K^{[n]}$. Let $\Delta$ be
the Young diagram of length $n+1$ such that the subvariety $X$
generically lies in the stratum $K^{[n]}_\Delta$ of the diagononal
stratification.

\punkt Assume first that $\Delta$ is not the diagram $\langle
1,1,\ldots,1\rangle$, so that the stratum $K^{[n]}_\Delta \subset
K^{[n]}$ is not the maximal one, and the dimension of the punctual
Hilbert scheme $F_\Delta$ is positive. Let the dense open subset
$\wh{U} \subset X \cap K^{[n]}_\Delta$ and the \'etale covering
$\wt{U} \to \wh{U}$ be as in \ref{wtU}. The map $\wt{U} \to \wh{U}
\subset K^{[n]}$ factors through the quotient map
$\eta:\wt{K}^{[n]}_\Delta \to K^{[n]}_\Delta$ by means of a map
$$
\rho = \rho_F \times \rho_T:\wt{U} \to F_\Delta \times \left(T^l_0
\setminus \Diag\right) = \wt{K}^{[n]}_\Delta,
$$
and by Proposition~\ref{point} the map $\rho_F:\wt{U} \to F_\Delta$
factors through a projection to a single point $o \in F_\Delta$. 

\punkt Let $a \in F_\Delta$ be any other point, and consider the subvariety
$$
U_a = \eta\left(\left(a \times
\rho_T\right)\left(\wt{U}\right)\right) \subset K^{[n]}_\Delta.  
$$
We claim the following. 

\begin{lemma}\label{closure}
The closure $\overline{U}_a \subset K^{[n]}$ is a
trianalytic subvariety in the Kummer variety $K^{[n]}$. 
\end{lemma}

\proof We first prove that the subset $\overline{U}_a \subset
K^{[n]}$ is analytic. Indeed, as established in the proof of
Proposition~\ref{point}, the closure $\overline{U} \subset T^l_0$ of
the subset $\rho_T(\wt{U}) \subset T^l_0 \setminus \Diag \subset
T^l_0$ is an analytic subvariety. Thus the closure in
$\overline{K}^{[n]}_\Delta = F_\Delta \times T^l_0$ of the subset
$$
\left(a \times \rho_T\right)\left(\wt{U}\right) = \{a\} \times
\rho_T\left(\wt{U}\right) \subset F_\Delta \times \left(T^l_0
\setminus \Diag\right),
$$
being equal to $\{a\} \times \overline{U} \subset F_\Delta \times
T^l_0$, is also analytic. Therefore $\overline{U}_a \subset K^{[n]}$
is analytic by Corollary~\ref{image} of Lemma~\ref{proper}.

Now, for every cohomology class $\alpha \in H^\idot(K^{[n]},\C)$
consider the Poincare pairing $\langle \alpha,
[\overline{U}_a]\rangle$ between $\alpha$ and the fundamental class
$[\overline{U}_a] \in H^\idot(K^{[n]},\C)$ of the analytic
subvariety $\overline{U}_a \subset K^{[n]}$. This pairing coincides
with the integral of a certain $C^\infty$ differential form
$\omega_\alpha$ over the dense open subset $U_a \subset
\overline{U}_a$. The map $\eta \circ (a \times \rho_T):\wt{U} \to
U_a$ is finite, say, of degree $m_a$, and we have
$$
m_a\left\langle \alpha, \left[\overline{U}_a\right]\right\rangle =
m_a\int_{U_a} \omega_\alpha = \int_{\wt{U}}\left(\eta \circ
\left(a \times \rho_T\right)\right)^*\omega_\alpha. 
$$
The right-hand side obviously depends continuously on the point $a
\in F_\Delta$. Since the Poincare pairing is non-degenerate, this
implies that the function
$$
a \mapsto m_a\left[\overline{U}_a\right] \in H^\idot(K^{[n]},\C)
$$
is continuous on the punctual Hilbert scheme $F_\Delta$. Since it
takes values inside the integral lattice $H^\idot(K^{[n]},{\Bbb Z})
\subset H^\idot(K^{[n]},\C)$, it must be locally constant. Moreover,
the punctual Hilbert scheme $F_\Delta$ is connected and irreducible
by \cite{Irr}. Therefore this function is constant on
$F_\Delta$. Thus 
$$
[\overline{U}_a] = \frac{m_o}{m_a}[\overline{U}_u] =
\frac{m_o}{m_a}[X] \in H^\idot(K^{[n]},\C),
$$
and the subvariety $\overline{U}_a \subset K^{[n]}$ is trianalytic
by the Trianalyticity Criterion of \cite{VerTrian} (Theorem
\ref{_triana_subva_SU(2)_Theorem_}).  
\endproof

\punkt Since by assumption $\dim F_\Delta > 0$, this lemma implies
that there exist a continuum of trianalytic subvarieties in the
Kummer variety $K^{[n]}$. They all have the same fundamental
class. Moreover, by our assumption on the maximality of the
dimension $\dim X$ and by Proposition~\ref{rigid_Proposition_}, all
these subvarieties are rigid. But the Douady moduli space of
subvarieties of the manifold $K^{[n]}$ with fixed fundamental class
is compact (\cite{_Fujiki_}, \cite{Lieber}). In particular, there
exists at most a finite number of rigid subvarieties, which is a
contradiction.

\rem The family of subvarieties $U_a \subset K^{[n]}_\Delta$
obviously depends holomorphically on the point $a \in F_\Delta$. A
finer argument should produce directly a structure of a holomorphic
family on the set of subvarieties $\overline{U}_a \subset K^{[n]}$,
thus avoiding references to the Compactness Theorem and to the
cardinality argument. However, it is not immediately clear how
to prove that taking closure preserves, at least generically, the
holomorphic dependence on the parameter $a$.

\punkt We have proved that the trianalytic subvariety $X \subset
K^{[n]}$ must necessarily lie generically in the maximal stratum
$$
K^{[n]}_{\langle 1,\ldots,1 \rangle} \subset K^{[n]}
$$ 
in the stratification by diagonals of the Kummer variety
$K^{[n]}$. To simplify notation, we will write $K^{[n]}_{gen}$ instead
of $K^{[n]}_{\langle 1,\ldots,1\rangle}$. The stratum $K^{[n]}_{gen}$ is
open and dense in $K^{[n]}$, and the map $\pi:K^{[n]} \to K^{(n)}$
is one-to-one on $K^{[n]}_{gen} \subset K^{[n]}$.

Let $U \subset X \cap K^{[n]}_{gen}$ be a smooth open dense subset in
the subvariety $X$. Let 
$$
\wt{U} = \eta^{-1}(U) \subset \wt{K}^{[n]}_{gen} = T^{n+1}_0 \setminus
\Diag
$$ 
be the preimage of the subset $U \subset K^{[n]}_{gen}$ under the Galois
covering $\eta:\wt{K}^{[n]}_{gen} \to K^{[n]}_{gen}$, and let 
$$
\overline{U} \subset T^{n+1}_0
$$
be its closure in the complex torus $T^{n+1}_0$. By the same
argument as in the proof of Lemma~\ref{closure}, the subset
$\overline{U} \subset T^{n+1}_0$ is a trianalytic
subvariety. Therefore by Lemma~\ref{generic} it is an image of a
union of hyperplanes in the universal covering $\C^{2n} \to
T^{n+1}_0$.

\punkt We finish the proof by repeating the argument used in the
case of a non-generic stratum $K^{[n]}_\Delta$. Namely, if $\dim X <
2n$, applying linear translation to this union of hyperplanes gives
a continuous family $U_a \subset \wt{K}^{[n]}_{gen}$ of different
subvarieties in $\wt{K}^{[n]}_{gen}$. This family contains $\wt{U}$ as
$U_a$ for some value of the parameter $a$. As in the non-generic
case, for every $a$ the closure $\overline{U}_a \subset K^{[n]}$ is
an analytic subvariety with the same fundamental class as $X \subset
K^{[n]}$. These subvarieties are all rigid, which,
is impossible since the Douady space is compact. 
Therefore we must have $\dim X = 2n$, or, in other
words, $X = K^{[n]}$. This finishes the proof of Theorem~\ref{main}.
\endproof

\end{document}